\documentclass[reqno,12pt]{article}
\usepackage{amssymb,amsthm,amsmath,amsfonts}
\usepackage{epsfig}

\theoremstyle{plain}
\begingroup

\newtheorem{theorem}{Theorem}[section]

\endgroup
\theoremstyle{definition}
\newtheorem{defn}{Definition}[section]

\newtheorem{example}[defn]{Example}

\theoremstyle{remark}

\numberwithin{equation}{section}
\numberwithin{figure}{section}

\DeclareMathOperator*{\res}{\mathrm{Res}}

\def\I{\mathrm{i}}

\def\D{{\mathbb D}}

\def\C{{\mathbb C}}
\def\P{{\mathbb P}}
\def\Z{{\mathbb Z}}

\begin{document}

\title{Line bundles defined by the Schwarz function
}

\author{
Bj\"orn Gustafsson\textsuperscript{1},
Mihai Putinar\textsuperscript{2}\\
}


\maketitle



\begin{center}
{\it To our friend Dima Khavinson, with best wishes}
\end{center}

\begin{abstract}
Cauchy and exponential transforms are characterized, and constructed, as canonical holomorphic
sections of certain line bundles on the Riemann sphere defined in terms of the Schwarz function.
A well known natural connection between Schwarz reflection and line bundles 
defined on the Schottky double of a planar domain is briefly discussed in the same context.
\end{abstract}



 \footnotetext[1]
{Department of Mathematics, KTH, 100 44, Stockholm, Sweden.\\
Email: \tt{gbjorn@kth.se}}
\footnotetext[2]
{Department of Mathematics, UCSB, Santa Barbara, CA 93106-3080, USA and School of Mathematics, Statistics and Physics,
Newcastle University, Newcastle upon Tyne, NE1 7RU, UK.\\
Email: \tt{mputinar@math.ucsb.edu, mihai.putinar@ncl.ac.uk}}



\section{Introduction}

The interpretation of function theoretic results in one complex variable as geometric statements on naturally associated Riemann surfaces is
often unveiling hidden and deeper structures. Schwarz reflection and its impersonation as Schwarz function of an analytic arc are not excluded from being enriched by this
higher, but very classical otherwise, perspective. It is quite surprising that global geometric interpretations of Schwarz function continue to be discovered
nowadays. The aim of the present
 note is to briefly recount two such instances.

In the first part we highlight a couple of discoveries made in \cite{Gustafsson-Putinar-2017}, viz. characterization of Cauchy and exponential transforms as unique holomorphic
sections of certain line bundles defined in terms of the Schwarz function. In contrast to the exposition in \cite{Gustafsson-Putinar-2017} we choose below a constructive approach,
starting from the line bundle in question and then finding its holomorphic sections by general methods, based on Cauchy integrals. In a different language, what we do is solving
some simple Riemann-Hilbert problems.

Despite our results being technically simple we find it striking that such a powerful tool as the exponential transform can be obtained, and characterized, 
just as a canonical section of an easily defined line bundle.
We recall that the exponential transform originally arose in operator theory via certain determinantal and trace formulas 
\cite{Pincus-1968}, \cite{Helton-Howe-1975}. See \cite{Martin-Putinar-1989}, \cite{Xia-2015} for overviews. 
It later became an important tool in the theory of quadrature domains  \cite{Putinar-1994},   \cite{Putinar-1995}, \cite{Gustafsson-Shapiro-2005},
it has been used (implicitly) for studying boundaries analytic varieties
 \cite{Harvey-Lawson-1975}, \cite{Alexander-Wermer-1998}, (explicitly) for proving regularity of free boundaries
\cite{Gustafsson-Putinar-1998}, for efficient reconstruction of domains from moment data \cite{Putinar-1996},  \cite{Gustafsson-He-Milanfar-Putinar-2000},
and for related Pad\'{e} approximation \cite{Putinar-2002}, to mention just a few instances. 

The Schwarz function \cite{Davis-1974}, \cite{Shapiro-1992} of an analytic curve is a more elementary object, 
which nevertheless has turned out to be highly useful in a variety of contexts in pure and applied mathematics. 
Among recent developments may be mentioned applications to gravitational lensing \cite{Fassnacht-Keeton-Khavinson-2009}, \cite{Khavinson-Neumann-2006}
and appearance of iterations of Schwarz functions in the study of the topology of quadrature domains \cite{Lee-Makarov-2016}.

The second part of the note is devoted to the relation between the Schwarz function of an analytic boundary in complex plane and various line bundles living on the Schottky double 
of the respective inner domain. The classical by now Riemann surface characterization of domains carrying various quadrature formulas for analytic functions is derived in these terms.


\section{Line bundles in general}\label{sec:line bundles}

Recall \cite{Gunning-1966}, \cite{Forster-1981}, \cite{Griffiths-Harris-1978} that a line bundle $\lambda$ on a compact Riemann surface $M$ is an element in the cohomology group
$H^1(M,\mathcal{O}^*)$, where $\mathcal{O}^*$ denotes the multiplicative sheaf of germs on nonvanishing
holomorphic functions. Using \v{C}ech cohomology this means that $\lambda$ is represented, in terms of
a sufficiently fine open covering $\{U_\alpha\}$ of $M$, by a collection of transition functions $\lambda_{\alpha\beta}\in\mathcal{O}^*(U_{\alpha}\cap U_{\beta})$,
one for each pair $\alpha,\beta$ with $U_{\alpha}\cap U_{\beta}\ne \emptyset$. These are subject to the multiplicative cocycle
conditions $\lambda_{\alpha\beta}\lambda_{\beta\alpha}=1$ and $\lambda_{\alpha\beta}\lambda_{\beta\gamma}\lambda_{\gamma\alpha}=1$
in $U_{\alpha}\cap U_{\beta}\cap U_{\gamma}$, whenever that intersection is nonempty.

There is also a notion of equivalence between two such collections $\lambda_{\alpha\beta}$ and $\lambda'_{\alpha'\beta'}$, reflecting at the level of cochains the same underlying line bundle $\lambda$, but we will 
actually be more interested in representing cocycles themselves rather than their equivalence classes in $H^1(M,\mathcal{O}^*)$.

In the present note we make some pertinent observations concerning the absolutely simplest case, namely the following situation:
$$
M=\P= \C\cup\{\infty\}=\text{the Riemann sphere},
$$
$$
\Gamma=\text{a closed analytic curve in} \,\,\C,
$$ 
which hence divides $\P$ into three disjoint pieces:
$$
\P=\Omega \cup \Gamma \cup \Omega^e.
$$ 
Here $\Omega$ is the bounded domain, hence $\infty\in\Omega^e$. We may assume, for later convenience,  that $0\in\Omega$.
Taking orientations into account we have
$$
\Gamma =\partial \Omega=-\partial \Omega^e.
$$

Let $U_1$ be a neighborhood of $\Omega\cup \Gamma$, $U_2$ a neighborhood of $\Omega^e\cup \Gamma$, and let $\lambda_{12}$ be analytic
and nonvanishing in $U_1\cap U_2$. In other words, $\lambda_{12}\in\mathcal{O}^*(U_1\cap U_2)$. Then $\lambda_{12}$ defines a line bundle
$\lambda\in H^1(\P,\mathcal{O}^*)$ according to the above scheme.

A meromorphic section of $\lambda$ is represented by a pair $f=(f_1,f_2)$ of meromorphic functions in $U_1$ and $U_2$, respectively, 
satisfying
\begin{equation}\label{flambdaf}
f_1=\lambda_{12} f_2 \quad \text{in}\,\, U_1\cap U_2.
\end{equation}
In situations like this we allow ourselves to shrink the initial open sets $U_1$ and $U_2$ if necessary. The {\it Chern class} of $\lambda$ is 
$$
c(\lambda)=\frac{1}{2\pi \I} \oint_\Gamma  d\log \lambda_{12},
$$ 
and this integer equals the number of zeros minus the number of poles for any meromorphic section  (not identically zero) of $\lambda$.
Of course, multiplicities should be taken into account, and in $U_1\cap U_2$ one counts for only one of the representatives $f_1$
and $f_2$. 

If $c(\lambda)<0$ there are no holomorphic sections of $\lambda$, while if $c(\lambda)\geq 0$ the vector space of holomorphic sections
has dimension $c(\lambda)+1$. In the present case, with the Riemann sphere, these are elementary statements, whereas corresponding 
(and more complicated) statements for higher genus Riemann surfaces require some global balance checking, well synthesized by the Riemann-Roch theorem. To make everything even more simple,
assume now that $c(\lambda)=0$. Then there is a one-dimensional space of holomorphic sections, and there will for example be a unique
holomorphic section of $\lambda$ which takes the value one at infinity (thus $f_2(\infty)=1$). How do we find this section $f=(f_1,f_2)$?

The construction is indeed quite easy (and of course well-known): write (\ref{flambdaf}) additively by taking the logarithm, as
\begin{equation}\label{logflambdaf}
\log f_1 =\log \lambda_{12} + \log f_2  \quad \text{in}\,\, U_1\cap U_2.
\end{equation}
Here $\log \lambda_{12}$ has a single-valued branch because $c(\lambda)=0$, and since $\Omega$ and $\Omega^e$ are simply connected
there are no problems with the other logarithms either. We may choose branches so that the requirement $f_2(\infty)=1$ means that
$\log f_2(\infty)=0$. Then, given $\lambda$, (\ref{logflambdaf}) is immediately solved for $f$ by using the Cauchy integral. Indeed, restricting 
(\ref{logflambdaf}) to $\Gamma$ it says that the two pieces of $\log f$ make the jump $\log \lambda_{12}$ across $\Gamma$, hence the solution
is obtained as the Cauchy integral:
\begin{equation}\label{Cauchylogflambdaf} 
\frac{1}{2\pi \I} \oint_\Gamma  \frac{\log \lambda_{12}(\zeta)d\zeta}{\zeta-z}=
\begin{cases}
\log f_1 (z), \quad z\in \Omega,\\
\log f_2 (z), \quad z\in \Omega^e.
\end{cases}
\end{equation}
It is easy to see that the two analytic functions $\log f_j$ extend analytically across $\Gamma$, hence the problem of constructing $f$ is solved. 

Corresponding problems in the cases $c({\lambda})\ne 0$ can likewise be solved, by small modifications of the above, in essence by taking
advantage of the equivalence between line bundles and divisor classes. 
One may for example rewrite (\ref{flambdaf}) as
\begin{equation}\label{adjustment}
f_1 (z) = \lambda_{12}(z)(z-a)^{-c(\lambda)}\cdot f_2(z) (z-a)^{c(\lambda)},
\end{equation}
for some fixed point $a\in\Omega$,
so that $\log (\lambda_{12}(z)(z-a)^{-c(\lambda)})$ is now single-valued on $\Gamma$. Then the Cauchy
integral corresponding to (\ref{Cauchylogflambdaf}) makes sense and gives
\begin{equation}\label{Cauchylogflambdafadjusted} 
\frac{1}{2\pi \I} \oint_\Gamma  \frac{\log( \lambda_{12}(\zeta) (\zeta-a)^{-c(\lambda)})d\zeta}{\zeta-z}=
\begin{cases}
\log f_1 (z), \quad& z\in \Omega,\\
\log (f_2 (z)(z-a)^{c(\lambda)}), \quad & z\in \Omega^e.
\end{cases}
\end{equation}


\section{The Cauchy transform}\label{sec:Cauchy}

Next we apply the procedure in Section~\ref{sec:line bundles} to a special case. We keep all notations and assumptions,
and recall that an analytic curve $\Gamma$ as above has a natural analytic function associated to it. This is the {\it Schwarz function}
$S(z)$ \cite{Davis-1974}, \cite{Shapiro-1992}, uniquely defined as the analytic function in a neighborhood of $\Gamma$ subject to the constraint 
\begin{equation}\label{S}
S(z)=\bar{z}, \quad z\in\Gamma.
\end{equation}
This means that $z\mapsto \overline{S(z)}$ is the anti-conformal reflection in $\Gamma$.

Consider the line bundle $\lambda$ defined by
\begin{equation}\label{lambdaeS}
\lambda_{12}(z)= e^{S(z)}
\end{equation}
in some neighborhood of $\Gamma$. It is immediately verified that $c(\lambda)=0$. So what is then the unique holomorphic section that takes the value one
at infinity?

The answer is directly clear from (\ref{Cauchylogflambdaf}). The Cauchy integral there becomes
$$
\frac{1}{2\pi \I} \oint_\Gamma  \frac{S(\zeta)d\zeta}{\zeta-z}=\frac{1}{2\pi \I} \oint_\Gamma  \frac{\bar{\zeta}d\zeta}{\zeta-z},
$$
and for $z\in \Omega^e$ this equals (minus) the Cauchy transform $C_\Omega(z)$ of $\Omega$ (i.e., of the measure which has density one in $\Omega$ and
vanishes elsewhere):
$$
\log f_2(z)=-\frac{1}{2\pi \I} \int_\Omega  \frac{d\zeta\wedge d\bar{\zeta}}{\zeta-z}
=-C_\Omega(z), \quad z\in\Omega^e.
$$
Similarly one may view the piece, $\log f_1(z)$, in $\Omega$ as a renormalized version of the Cauchy transform of $\Omega^e$.
We may define this exterior Cauchy transform simply as the corresponding boundary integral above:
$$
\log f_1(z) =\frac{1}{2\pi \I} \oint_\Gamma  \frac{\bar{\zeta}d\zeta}{\zeta-z}= C_{\Omega^e}(z), \quad z\in\Omega.
$$
In summary, the exponential of the Cauchy transforms associated to $\Gamma$ defines the unique holomorphic section of 
$\lambda=e^S$ taking the value one at infinity.

Expanding the germs of $\log f_j$ at $z=0$ and $z=\infty$ in powers of $z$ gives
$$
\log f_1(z) =\sum_{k<0} \frac{M_k}{z^{k+1}},
$$
$$
\log f_2(z) =-\sum_{k\geq 0} \frac{M_k}{z^{k+1}},
$$
where the $M_k$ are the {\it harmonic moments} of $\Omega$,  defined for arbitrary $k\in \Z$ by
$$
M_k=\frac{1}{2\pi \I} \oint_\Gamma  {z^k \bar{z}dz}.
$$
Thus the relation (\ref{logflambdaf}), written as
$$
 \log \lambda_{12} =\log f_1-  \log f_2  \quad \text{in}\,\, U_1\cap U_2,
$$
can in the present situation be identified with the well-known expansion
$$
S(z)=\sum_{k\in \Z} \frac{M_k}{z^{k+1}}. 
$$
of the Schwarz function. Note however that the full series here need not converge anywhere, it is just a formal sum of two germs of analytic functions.


\section{The exponential transform}\label{sec:exponential}

Finally we apply the construction in Section~\ref{sec:line bundles}  to 
line bundles defined in terms of the Schwarz function $S(z)$ and a parameter $w\in\C\setminus\Gamma$.
We define these bundles $\lambda=\lambda_w$  by the transition functions
\begin{equation}\label{lambda12zw}
\lambda_{12}(z)=\frac{1}{S(z)-\bar{w}}.
\end{equation}
The reciprocals are even simpler:
$$
\lambda_{21}(z)={S(z)-\bar{w}}.
$$

We first choose $w\in\Omega^e$.
Then one finds $c({\lambda})=0$, so there is a unique holomorphic section taking the value one at infinity.
To find it we again apply (\ref{Cauchylogflambdaf}). The exterior part ($z\in \Omega^e$) of the Cauchy integral there becomes
\begin{equation}\label{logf2}
\log f_2 (z)= -\frac{1}{2\pi \I} \oint_\Gamma  \frac{\log (S(\zeta)-\bar{w})d\zeta}{\zeta-z}=-\frac{1}{2\pi \I} \oint_\Gamma  \frac{\log (\bar{\zeta}-\bar{w})d\zeta}{\zeta-z}
\end{equation}
\begin{equation}\label{Czw}
=\frac{1}{2\pi \I} \int_\Omega  \frac{d{\zeta}\wedge d\bar{\zeta}}{(\zeta-z)(\bar{\zeta}-\bar{w})}
=C_\Omega(z,w).
\end{equation}
The last equality sign here defines the {\it double Cauchy transform} $C_\Omega(z,w)$ of $\Omega$.
Its exponential is the {\it exponential transform} of $\Omega$:
\begin{equation}\label{EC}
E_\Omega(z,w)=\exp C_\Omega(z,w).
\end{equation}
We let the functions $C_\Omega(z,w)$ and $E_\Omega(z,w)$ be defined by (\ref{Czw}), (\ref{EC}) in all $\C\times\C$.
Clearly they are hermitian symmetric in $z$ and $w$. In addition,  $-C_\Omega(z,w)$ and $1/E_\Omega(z,w)$ can be seen to be
positive definite in $\Omega^e$ (see \cite{Gustafsson-Putinar-2017} for details, and for further positivity properties).

Thus we infer
\begin{equation}\label{f2E}
f_2(z) ={E_\Omega(z,w)}, \quad z, w\in \Omega^e.
\end{equation}
Note that $f_2$ depends on $w$ via $\lambda=\lambda_w$. As a function of $z$ it remains analytic in $U_2\supset \Omega^e \cup \Gamma$,
and, because of the hermitian symmetry, $E_\Omega(z,w)$ has a corresponding anti-analytic extension also in the variable $w$. We shall denote
the full analytic/anti-analytic extension of $E_\Omega(z,w)$, in both variables,  by $F(z,w)$,  to distinguish it from $E_\Omega(z,w)$ itself, which is  defined everywhere 
by (\ref{Czw}), (\ref{EC}). 
Consequently
$$
f_2(z)=F(z,w), \quad z\in U_2, \,\, w\in\Omega^e.
$$

For $\log f_1(z)$, where then $z\in\Omega$, the Cauchy integral in (\ref{logf2}) gets an additional contribution at $\zeta=z$. The result is
$$
\log f_1(z) =C_\Omega(z,w)-\log (\bar{z}-\bar{w}), \quad z\in \Omega,\,\, w\in \Omega^e,
$$ 
$$
f_1(z) =\frac{E_\Omega(z,w)}{\bar{z}-\bar{w}}=G(z,w), \quad z\in \Omega,\,\,  w\in \Omega^e.
$$
Here the last equality is simply the definition of the function $G(z,w)$, which is analytic/anti-analytic in $\Omega\times \Omega^e$.

In summary, $(f_1,f_2)=(G(\cdot,w), F(\cdot,w))$ defines the unique normalized holomorphic section of the bundle $\lambda=\lambda_w$
when $w\in \Omega^e$. Explicitly the transition relation is
\begin{equation}\label{GSE}
G(z,w)(S(z)-\bar{w})=F(z,w), \quad z\in U_1\cap U_2, \,\,w\in \Omega^e.  
\end{equation}

When $w\in \Omega$, the bundle defined by (\ref{lambda12zw}) has Chern class
$$
c(\lambda)=1.
$$
This means that $\lambda$ is equivalent to the {\it hyperplane section line bundle}  (see \cite{Griffiths-Harris-1978} for the terminology), for which the space 
of holomorphic sections has dimension two,  and there is  a unique such section which vanishes with leading term $1/z$ at infinity. 
From this we obtain two new analytic pieces ($G^*$ and $H$ below) of the exponential transform, by using the adjusted
Cauchy integral (\ref{Cauchylogflambdafadjusted}) with $a=w$. Considering first $\log f_2(z)$, hence having $z\in\Omega^e$, we get
$$
\log f_2 (z) (z-w)= -\frac{1}{2\pi \I} \oint_\Gamma  \frac{\log (S(\zeta)-\bar{w})(\zeta-w)d\zeta}{\zeta-z}
$$
$$
=-\frac{1}{2\pi \I} \oint_\Gamma  \frac{\log (\bar{\zeta}-\bar{w}) (\zeta-w)d\zeta}{\zeta-z}.
$$
($\log fg$ is to be interpreted as $\log (fg)$, in general.) This again equals $C_\Omega(z,w)$
since the logarithmic singularities at $\zeta=w$ do not give any distributional contributions (or, the line integral around
$|\zeta-w|=\varepsilon$ tends to zero as $\varepsilon\to 0$). Hence
$$
f_2(z)= \frac{E_\Omega(z,w)}{z-w}=- \overline{G(w,z)} =-G^*(z,w),
$$
where $G^*(z,w)=\overline{G(w,z)}$, by definition. It is easy to see that $f_2(z)$ behaves as $1/z$ at infinity, which is
the normalization we had set up. 

For $\log f_1(z)$, hence with $z\in\Omega$, we get as before an additional term in the right member 
because of the singularity at $\zeta=z$:
$$
\log f_1 (z)= -\frac{1}{2\pi \I} \oint_\Gamma  \frac{\log (\bar{\zeta}-\bar{w})(\zeta-w)d\zeta}{\zeta-z}
$$
$$
= C_\Omega(z,w)-\log (z-w)(\bar{z}-\bar{w}).
$$
Thus
$$
f_1(z)=\frac{E_\Omega (z,w)}{|z-w|^2}=H(z,w), \quad z,w \in \Omega,
$$
where the last equality defines the {\it interior exponential transform}  $H(z,w)$.

So when $w\in \Omega$, the line bundle defined by (\ref{lambda12zw}) has a uniques holomorphic section behaving as $1/z$
at infinity, and this is given by $(f_1,f_2)=(H(\cdot, w), -G^*(\cdot, w))$.  The explicit transition formula  becomes
\begin{equation}\label{HSG}
H(z,w)(S(z)-\bar{w})=-G^*(z,w), \quad z\in U_1\cap U_2, \,\,w\in \Omega.  
\end{equation}

\begin{example}
With $\Omega=\D$ we have
$$
S(z)-\bar{w}= \frac{1}{z}-\bar{w},
$$
and the relations (\ref{GSE}) and (\ref{HSG}) become
\begin{align*}
(-\frac{1}{\bar{w}})\cdot (\frac{1}{z}-\bar{w})&=1-\frac{1}{z\bar{w}},\\
\frac{1}{1-z\bar{w}}\cdot(\frac{1}{z}-\bar{w})&= \frac{1}{z},
\end{align*}
thereby also exhibiting the functions $F$, $G$, $G^*$, $H$ explicitly in this case.
\end{example}

We summarize the observations accumulated in the previous section as follows.

\begin{theorem}
Let $\Omega\subset\C$ be a bounded simply connected domain for which $\Gamma=\partial\Omega$
is a smooth analytic curve, let $S(z)$ be the Schwarz function of $\Gamma$ and let, for each $w\in \C\setminus \Gamma$,
$\lambda=\lambda_w$ be the line bundle on the Riemann sphere defined by (\ref{lambda12zw}).

The above defined four analytic constituents $F(z,w)$, $G(z,w)$, $G^*(z,w)$, $H(z,w)$ 
of the exponential transform $E_\Omega (z,w)$ of $\Omega$ (see (\ref{Czw}), (\ref{EC})) 
paste together, in each variable separately, to canonically defined holomorphic sections of $\lambda$, as follows.

\begin{itemize}

\item For each fixed $w\in \Omega^e$, the pair $(G(\cdot,w), F(\cdot, w))$ represents the unique holomorphic section of $\lambda_w$ which takes the value one at infinity.  
Explicitly this means that (\ref{GSE}) holds and  characterizes the pair uniquely under $F(\infty,w)=1$. 

\item For each fixed $w\in \Omega$, the pair $(H(\cdot,w), -G^*(\cdot, w))$ represents the unique holomorphic section of $\lambda_w$ which vanishes at infinity with leading
term $1/z$.  
Explicitly this means that (\ref{HSG}) holds and characterizes the pair uniquely under $-G^*(z,w)=z^{-1}+\mathcal{O}(z^{-2})$, $z\to\infty$. 

\end{itemize}

As a limiting case of the above, when $w\to \infty$, we infer that the exponentials of the Cauchy transforms of $\Omega$ and $\Omega^e$ paste together to represent
the unique holomorphic section of the line bundle defined by (\ref{lambdaeS}) taking the value one at infinity (see Section~\ref{sec:Cauchy} for details).

\end{theorem}


\section{The Schwarz function and the Schottky double}

This section is aimed at reviewing from a global geometric perspective the interplay between Schwarz's function and line bundles on Schottky doubles of planar domains.

\subsection{Line bundles on the Schottky double}

If $\Omega\subset \C$ is a, possibly multiply connected, domain with analytic boundary, then its {\it Schottky double} is the compact symmetric Riemann surface $\hat\Omega$
obtained by completing $\Omega$ with a backside $\tilde{\Omega}$, a copy of $\Omega$ provided with the opposite conformal structure, the two
surfaces being glued  together along their common boundary $\Gamma=\partial\Omega$. 
Denoting by $\tilde{z}$ the point on $\tilde{\Omega}$ corresponding to $z\in\Omega$, a natural holomorphic atlas for $\hat{\Omega}$ is $\{(U_j,\phi_j):j=1,2\}$, where the two
coordinate maps $\phi_j:U_j\to\C$ are initially defined by
$$
\phi_1:z\mapsto z \quad (z\in\Omega), 
$$
$$
\phi_2: \tilde{z}\mapsto \bar{z}\quad  (\tilde{z}\in\tilde{\Omega}),
$$
and are then,  using the analyticity of $\Gamma$, extended analytically across $\Gamma$ in $\hat{\Omega}$ to some open sets
$U_1\supset \Omega\cup\Gamma$ and $U_2\supset \Omega^e \cup \Gamma$. Here $U_1$ is to be kept so small that $\phi_1$ is
still univalent in $U_1$. 
Now, the coordinate transition function between $\phi_1$ and $\phi_2$ is exactly the Schwarz function:
$$
S=\phi_2 \circ \phi_1^{-1}.
$$ 
Indeed, $S$ is defined and holomorphic in $\phi_1(U_1\cap U_2)\subset\C$, and (\ref{S}) holds.

From this it is clear that several line bundles can be defined via $S$. First of all we have the {\it canonical line bundle}, with the transition function
$$
\kappa_{12}=\frac{d\phi_2}{d\phi_1}.
$$
Thinking of this as being defined on $\phi_1(U_1\cap U_2)$, rather than on $U_1\cap U_2$, it is simply the derivative of
the Schwarz function:
$$
S'=\kappa_{12}\circ \phi_1^{-1}.
$$
Sections of $\kappa$ are differential forms (of type $(1,0)$).  The most handy way to represent such a section is to let the two pieces, $f_1$ and $f_2$ say, 
both be defined in $\Omega$ in such a way that the actual functions on $U_1$ and $U_2$ are $f_1\circ \phi_1$ and ${\rm conj}\circ f_2\circ {\rm conj} \circ \phi_2$,
where ${\rm conj}$ denotes complex conjugation. Then the matching condition becomes
$$
f_1(z)=S'(z)\overline{f_2(z)} \quad \text{on\,\,} \Gamma, 
$$ 
more intuitively written as
$$
f_1 (z)dz =\overline{f_2 (z)dz} \quad \text{along}\,\, \Gamma.
$$
   
Denoting by $\texttt{g}$ the number of inner components of $\Gamma$, i.e., the number of ``holes'' of $\Omega$,
the Chern class of $\kappa$ is 
$$
c(\kappa)= \frac{1}{2\pi i} \oint_\Gamma d\log S' (z)= \frac{1}{2\pi i} \oint_\Gamma d\log \frac{d\bar{z}}{dz}=2( \texttt{g}-1),
$$
where the term ``$-1$'' comes from the outer component of $\Gamma$. 
The result  is of course what is  expected from general theory,  since $\texttt{g}$ is also the genus of $\hat{\Omega}$.

To define more bundles, let $T(z)$ denote the unit tangent vector on $\Gamma$, holomorphically extended to a neighborhood of $\Gamma$.
This means that
\begin{equation}\label{ST}
S'(z)=\frac{1}{T(z)^2}.
\end{equation}
Thus we get the possibility of defining a unique square root of $S'(z)$, and so a bundle $\lambda$ with $\lambda^2={\kappa}$.
We just set
\begin{equation}\label{lambdaT}
\lambda_{12}=\frac{1}{T\circ \phi_1}.
\end{equation}

It is understood above that $T(z)$ is oriented 
according to the standard orientation of $\Gamma$ as the boundary of $\Omega$,
but actually (\ref{ST}) holds also if one changes sign
of  $T(z)$, arbitrarily on each component of $\Gamma$.  In this way one can get $2^{\texttt{g}}$ different bundles $\lambda$ satisfying $\lambda^2=\kappa$.
Actually there are more such bundles (precisely $2^{2\texttt{g}}$, see \cite{Hejhal-1972}), but most of them can only be defined by using a finer atlas than the one we have used.
Bundles $\lambda$ satisfying $\lambda^2=\kappa$ have been studied at length in several papers, for example \cite{Hawley-Schiffer-1966}, \cite{Hejhal-1972},
and in physics they are often referred to as {\it spin bundles} (see \cite{Atiyah-1971}, for example).

Both $\kappa$ and $\lambda$ have interesting canonical sections. For the double of a planar domain, and choosing the standard orientation for $T(z)$ in (\ref{lambdaT}),
there is a unique (after normalization) meromorphic section which has a simple pole at any given point $\tilde{a}\in\tilde{\Omega}$ and is elsewhere holomorphic.
The restriction of this section to $\Omega$ is the classical Szeg\"o kernel, reproducing the value at $a\in\Omega$ of any analytic function in the appropriate
Hilbert space  (Hardy, or Smirnov, space \cite{Duren-1970}).  A similar statement can be made for $\kappa$ and the Bergman kernel.
For some recent discussion of these classical matters, see \cite{Bell-2016}.

One can certainly take arbitrary integer powers $\lambda^m$ ($m\in\Z$) of $\lambda$, and the corresponding sections (differentials of order $m/2$) are most conveniently described, as above, 
by pairs $(f_1,f_2)$ of functions  in $\Omega$ satisfying the matching condition
\begin{equation}\label{f1f2T}
f_1(z)T(z)^m=\overline{f_2(z)} \quad \text{on\,\,} \Gamma.
\end{equation} 
Clearly, the Chern class of the bundle $\lambda^m$ is $c(\lambda^m)= m(\texttt{g}-1)$.


 \subsection{Sections of line bundles and quadrature domains}

Next we turn to some special situations in which functions related to the Schwarz function are meromorphic sections of  line bundles on the Schottky double.

The coordinate function $\phi_1$, which is holomorphic in  $U_1\subset \hat{\Omega}$, certainly cannot be extended to be  holomorphic in all $\hat{\Omega}$, 
but it may very well be extended to be meromorphic in $\hat{\Omega}$. This occurs if and
only if $\Omega$ is a {\it (classical) quadrature domain} (see \cite{Shapiro-1992},
\cite{Ebenfelt-Gustafsson-Khavinson-Putinar-2005}, \cite{Gustafsson-Shapiro-2005} for overviews),
or by other names, a {\it finitely determined domain} \cite{Putinar-1994} or and {\it algebraic domain} \cite{Varchenko-Etingof-1992}. 
An equivalent statement is that the Schwarz function 
extends to be meromorphic in all $\Omega$, in which case (\ref{S}) becomes an instance of (\ref{f1f2T}) with $m=0$.
Having $S(z)$ meromorphic in $\Omega$ immediately leads to a finite quadrature for analytic functions $f$ in $\Omega$, as noticed already in \cite{Davis-1974}:
$$
\frac{1}{2\pi i} \int_\Omega f(z)d\bar{z}dz =\frac{1}{2\pi i} \int_{\partial\Omega} f(z)S(z)dz= \sum_\Omega \res (f(z)S(z)dz).
$$

When $\phi_1$ extends meromorphically, so does $\phi_2$ and there will be an algebraic (polynomial) equation $Q(\phi_1, \phi_2)=0$ relating the two. Equivalently,
$Q(z,S(z))=0$ (identically). In particular, $Q(z,\bar{z})=0$ on $\partial\Omega$, i.e., boundaries of quadrature domains are algebraic, as was discovered in \cite{Aharonov-Shapiro-1976}.

We may at this point return to the exponential transform and notice that assuming $S(z)$ is meromorphic in $\Omega$ implies that the entire left hand side in (\ref{GSE}) is meromorphic
in $\Omega$, while the right hand side is holomorphic. Thus the two sides match together to define a meromorphic function in all $\P$, in other words a rational function.
The conclusion is that $F(z,w)$ is a rational function  in $z$ when $\Omega$ is a quadrature domain, and by hermitian symmetry it is also rational in $\bar{w}$. 
The more specific structure of this rational function turns out to be
\begin{equation}\label{FQPP}
F(z,w)= \frac{Q(z,\bar{w})}{P(z)\overline{P(w)}},
\end{equation}
where $Q$ is the same polynomial as above and $P(z)$ is a polynomial vanishing exactly  at the poles of $S(z)$.
See \cite{Putinar-1994}, \cite{Putinar-1996}, \cite{Gustafsson-Putinar-2017}, for example.

The above statements actually also go the other way around, in a quite strong form: if the germ of $E_\Omega(z,w)=F(z,w)$ at infinity (in both variables) is a rational 
function of the form (\ref{FQPP}) then $\Omega$ is a quadrature domain. (In Section~\ref{sec:exponential} it was assumed that $\Omega$ is simply connected, 
but this assumption is in no way essential.) 

Besides $\phi_1$ itself extending meromorphically to $\hat{\Omega}$ there are some other possibilities, for example that $d\phi_1$ extends meromorphically as a section of $\kappa$,
or that $\sqrt{d\phi_1}$ extends meromorphically as a section of $\lambda=\sqrt{\kappa}$. These give rise to other kinds
of quadrature domains. Efficient descriptions of these options can preferably be made in terms of (\ref{f1f2T}), as below.

\begin{itemize}

\item[(i)] The case of {\it classical quadrature domains} corresponds to $m=0$, $f_1(z)=z$, $f_2=S(z)$ (with $S(z)$ meromorphic in $\Omega$).

\item[(ii)] That $d\phi_1$ extends meromorphically is the same as saying that $S'(z)$ is meromorphic in $\Omega$, while then $S(z)$ itself may have logarithmic poles  and be additively
multivalued. This becomes the case $m=2$, $f_1(z)=1$, $f_2(z)=S'(z)= T(z)^{-2}$ in (\ref{f1f2T}) and gives rise to the finite quadrature
$$
\frac{1}{2\pi i} \int_\Omega f'(z)d\bar{z}dz =-\frac{1}{2\pi i} \int_{\partial\Omega} f(z)S'(z)dz =- \sum_\Omega \res (f(z)S'(z)dz).
$$
Such domains $\Omega$ are called {\it Abelian domains} in \cite{Varchenko-Etingof-1992}.

\item[(iii)] That $\sqrt{d\phi_1}$ extends, finally, is the same as saying that $\sqrt{S'(z)}=T(z)^{-1}$ is meromorphic in $\Omega$, which is represented by $m=1$, $f_1(z)=1$,
$f_2(z)= T(z)^{-1}$ in (\ref{f1f2T}). The resulting finite quadrature is
 $$
 \int_{\partial\Omega} f(z)|dz| =   \int_{\partial\Omega} f(z)\frac{dz}{T(z)} =2\pi i \sum_\Omega \res \frac{f(z)dz}{T(z)}.
$$
Such quadrature domains are discussed at length in \cite{Gustafsson-1987}, exactly from the point of view of half order differentials. 
(Corresponding studies for the cases (i) and (ii) were made in \cite{Gustafsson-1983}.)
Previous work include \cite{Shapiro-Ullemar-1981}, \cite{Avci-1977}. See again \cite{Gustafsson-Shapiro-2005} for a relatively recent overview.

\end{itemize}

In this context we may mention, finally, a finite analytic quadrature for polygons, discovered by Motzkin and Schoenberg and later refined by Philip Davis \cite{Davis-1964}.
For polygons the Schwarz function is piecewise linear, with changes of representatives at the corners. Thus, near the boundary, $S''(z)$ is zero except for Dirac contributions
at the corners. By partial integration on the boundary this gives a formula of the kind
$$
\frac{1}{2\pi i} \int_\Omega f''(z)d\bar{z}dz
=\frac{1}{2\pi i} \int_{\partial\Omega} f(z)S''(z)dz=\sum_j c_j f(a_j),
$$
where $a_j\in \partial\Omega$ are the corners.


\bibliography{bibliography_gbjorn.bib}

\end{document}